\documentclass[a4paper,reqno]{amsart}
\usepackage{amsmath}
\usepackage{amssymb}
\usepackage{hyperref}
\usepackage{url}
\usepackage{enumitem}
\usepackage{latexsym}
\usepackage{mathrsfs}
\usepackage[all]{xy}
\input{diagxy}

\hypersetup{
    pdfmenubar=false,
    pdffitwindow=true,
    pdfstartview=FitH,
    colorlinks=true,
    linkcolor=red,
    citecolor=green,
    urlcolor=cyan
}

\theoremstyle{plain}
\newtheorem{theorem}{\bf Theorem}[section]
\newtheorem{proposition}[theorem]{\bf Proposition}
\newtheorem{lemma}[theorem]{\bf Lemma}

\theoremstyle{definition}
\newtheorem{example}[theorem]{\bf Example}

\newtheorem{definition}[theorem]{\bf Definition}
\newtheorem{remark}[theorem]{\bf Remark}

\newcommand{\N}{\mathbb N}
\newcommand{\Z}{\mathbb Z}

\newcommand{\Q}{\mathbb Q}

\newcommand{\DP}{\negthinspace : \negthinspace}
\renewcommand{\time}{\negthinspace \times \negthinspace}

\newcommand{\red}{\text{\rm red}}

\newcommand{\id}{\text{\rm id}}
\newcommand{\free}{\mathcal{F}}
\newcommand{\block}{\mathcal{B}}
\newcommand{\st}{\,\mid\,}

\DeclareMathOperator{\coker}{Coker}
\DeclareMathOperator{\spec}{spec} \DeclareMathOperator{\supp}{supp}

\DeclareMathOperator{\End}{End}

\newcommand{\cdv}{\negthinspace \cdot_v \negthinspace}
\newcommand{\vspec}{\text{$v$-\rm{spec}}}
\newcommand{\ve}{\varepsilon}
\newcommand{\DivMon}{\mathcal{I}_v^*}

\renewcommand{\t}{\, | \,}

\renewcommand{\time}{\negthinspace \times \negthinspace}

\renewcommand{\to}{\rightarrow}

\numberwithin{equation}{section}

\begin{document}

\title[]
{On transfer homomorphisms of Krull monoids}

\author{Alfred Geroldinger and Florian Kainrath}

\address{University of Graz, NAWI Graz \\
Institute of Mathematics and Scientific Computing \\
Heinrichstra{\ss}e 36\\
8010 Graz, Austria}
\email{alfred.geroldinger@uni-graz.at,  florian.kainrath@uni-graz.at}
\urladdr{https://imsc.uni-graz.at/geroldinger}

\thanks{This work was supported by the Austrian Science Fund FWF, Project P33499-N}



\begin{abstract}
Every Krull monoid has a transfer homomorphism onto a monoid of zero-sum sequences over a subset of its class group. This transfer homomorphism is a crucial tool for studying the arithmetic of Krull monoids. In the present paper, we strengthen and refine this tool for Krull monoids with finitely generated class group.
\end{abstract}

\maketitle


\section{Introduction}

Transfer homomorphisms are a central tool in factorization theory. Since they are essentially surjective and allow to lift factorizations, they make it possible to study the arithmetic of monoids and domains as follows. Given an object $H$ of interest,   find a simpler object $B$, and a transfer homomorphism  $\theta \colon H \to B$; study the arithmetic of the simpler object $B$ and then  pull back arithmetical results from $B$ to the original object of interest $H$.  Oftentimes, the object $B$ is constructed with divisor theoretic or ideal theoretic methods and it has the same algebraic structure as $H$ but a simpler combinatorial structure. This holds true, in particular, for Krull monoids, weakly Krull monoids, and C-monoids, and for these classes of monoids the described strategy has turned out to be highly efficient. We refer to \cite{Ge-HK06a} for background on transfer homomorphisms.

In the present paper we focus on transfer homomorphisms of Krull monoids. Since a domain $D$ is a Krull domain if and only if its multiplicative monoid $D \setminus \{0\}$ of nonzero elements is a Krull monoid, we also have the domain case in mind when speaking of Krull monoids. Let $H$ be a Krull monoid with divisor theory $\varphi \colon H \to \mathcal F (P)$, class group $G$,  and let $G_P \subset G$ denote the set of classes containing prime divisors. Then the monoid $\mathcal B (G_P)$ of zero-sum sequences over $G_P$ is a Krull monoid again, and there is a canonical transfer homomorphism $\theta \colon H \to \mathcal B (G_P)$. It is an easy observation to see that the inclusion $\mathcal B (G_P) \hookrightarrow \mathcal F (G_P)$, where $\mathcal F (G_P)$ is the free abelian monoid with basis $G_P$, is a cofinal divisor homomorphism but in general not a divisor theory. This means that $\mathcal B (G_P)$ is not optimal, neither from the algebraic point of view nor for arithmetical investigations.

It turned out, first for special subsets of finite abelian groups, that it is possible to construct a new subset $G_P'$  and a divisor homomorphism $\theta' \colon H \to \mathcal B (G_P')$ such that $\mathcal B (G_P') \hookrightarrow  \mathcal F (G_P')$ is a divisor theory. Such refinement constructions are used in a variety of
 arithmetical investigations of $H $ (for a sample, see \cite[Theorem 3.17]{Sc04a}, \cite[Theorem 3.1]{Sc05d}, \cite[Chapter 6.7]{Ge-HK06a}, \cite[Lemma 4.5]{Pl-Sc20a}). In \cite{Sc10a}, Schmid proved that for every Krull monoid $H$ with torsion class group $G$ there is a  subset $G_0^* \subset G$, and a transfer homomorphism $\theta^* \colon H \to \mathcal B (G_0^*)$ such that $\mathcal B (G_0^*) \hookrightarrow \mathcal F (G_0^*)$ is a divisor theory (Proposition \ref{3.3}). This was done by transfinite induction over a combinatorial construction. In  this paper, we show the same result for Krull monoids with finitely generated class group (Theorem \ref{T:MainThm}). In contrast to the result by Schmid, we use ideal and divisor theoretic methods.

In Section \ref{2}, we gather the necessary background on Krull monoids. Our main result (Theorem \ref{T:MainThm}) is proved in Section \ref{3}. In the last section, we apply our results to transfer Krull monoids (Theorem \ref{4.3}) and to weakly Krull monoids (Theorem \ref{4.5}). Finally, we demonstrate the usefulness of our construction in an example of a Krull monoid with infinite class group (Example \ref{4.7}).

\smallskip
\section{Background on Krull monoids} \label{2}
\smallskip

\smallskip
\noindent
{\bf Arithmetic of monoids.} By a monoid, we mean a semigroup with identity element. If not stated otherwise, we use multiplicative notation. Let $H$ be a monoid with identity $1_H \in H$, and let $H^{\times}$ denote the group of invertible elements. We say that $H$ is
\begin{itemize}
\item {\it unit-cancellative} if $a, u \in H$ and $a=au$ or $a=ua$ implies $u \in H^{\times}$;

\item {\it cancellative} if $a, b, c \in H$ and $ab=ac$ or $ba=ca$ implies $b=c$.
\end{itemize}
A non-unit $a \in H$ is said to be an {\it atom} (or {\it irreducible}) if $b, c \in H$ and $a = bc$ implies that $b \in H^{\times}$ or $c \in H^{\times}$. We denote by $\mathcal A (H)$ the set of atoms of $H$ and we say that $H$ is atomic if every non-unit is a finite product of atoms. If $a = u_1 \cdot \ldots \cdot u_k \in H$ is a product of $k$ atoms, then $k$ is called a factorization length of $a$ and the set $\mathsf L_H (a) = \mathsf L (a) \subset \N$ is the {\it set of lengths} of $a$. It is convenient to set $\mathsf L (u)= \{0\}$ for all $u \in H^{\times}$. Then
\[
\mathcal L (H) = \{\mathsf L (a) \mid a \in H \}
\]
denotes the {\it system of sets of lengths} of $H$ and
\[
\Delta (H) = \bigcup_{L \in \mathcal L (H)} \Delta (L) \ \subset \N
\]
denotes the {\it set of distances} of $H$, where $\Delta (L)$ is the set of successive distances of $L$. For a set $P$, let $\mathcal F (P)$ denote the free abelian monoid with basis $P$. Then every $a \in \mathcal F (P)$ has a unique representation in the form
\[
a = \prod_{p \in P} p^{\mathsf v_p (a) } \,,
\]
where $\mathsf v_p \colon \mathcal F (P) \to \N_0$ is the $p$-adic valuation of $a$. We call $|a|= \sum_{p \in P} \mathsf v_p (a) \in \N_0$ the {\it length} of $a$ and $\supp (a) = \{ p \in P \mid \mathsf v_p (a) > 0 \} \subset P$ the {\it support} of $a$.

Now suppose that $H$ is a commutative unit-cancellative monoid. We say that $H$ is reduced if $H^{\times} = \{1_H\}$, and we denote by  $H_{\red} = \{a H^{\times} \mid a \in H\}$  the associated reduced monoid of $H$. The free abelian monoid $\mathsf Z (H) = \mathcal F ( \mathcal A (H_{\red}))$ is the factorization monoid of $H$, and we denote by $\pi \colon  \mathsf Z (H) \to H_{\red}$ the canonical epimorphism. For an element $a \in H$,
\[
\mathsf Z (a) = \pi^{-1} (aH^{\times}) \subset \mathsf Z (H) \quad \text{is the {\it set of factorizations} of $a$} \,.
\]
To define a distance function on $\mathsf Z (H)$, let
\[
z = u_1 \cdot \ldots \cdot u_kv_1 \cdot \ldots \cdot v_{\ell} \quad \text{and} \quad z' = u_1 \cdot \ldots \cdot u_kw_1 \cdot \ldots \cdot w_m
\]
be two factorizations, where $k, \ell, m \in \N_0$, all $u_r, v_s, w_t \in \mathcal A (H_{\red})$ such that $v_s \ne w_t$ for all $s \in [1, \ell]$ and $t \in [1, m]$. Then $\mathsf d (z, z') = \max \{\ell, m\} \in \N_0$ is the distance between $z$ and $z'$. Let $a \in H$ and $M \in \N_0$. A sequence $z_0, \ldots, z_k \in \mathsf Z (a)$ is called an $M$-chain of factorizations if $\mathsf d (z_{i-1}, z_i) \le M$ for all $i \in [1,k]$. The catenary degree $\mathsf c (a)$ is the smallest $M  \in \N_0 \cup \{\infty\}$ such that any two factorizations $z, z' \in \mathsf Z (a)$ can be concatenated by an $M$-chain. Then
\[
\mathsf c (H) = \sup \{\mathsf c (a) \mid a \in H \} \in \N_0 \cup \{\infty\}
\]
is the {\it catenary degree} of $H$. For $u \in \mathcal A (H_{\red})$, the {\it local tame degree} $\mathsf t (H, u)$ is the smallest $M \in \N_0 \cup \{\infty\}$ having the following property: if $a \in H$ with $\mathsf Z (a) \cap u \mathsf Z (H) \ne \emptyset$ and $z \in \mathsf Z (a)$, then there is $z' \in \mathsf Z (a) \cap u \mathsf Z (H)$ such that $\mathsf d (z, z') \le M$. The monoid $H$ is called
\begin{itemize}
\item {\it locally tame} if $\mathsf t (H, u) < \infty$ for all $u \in \mathcal A (H_{\red})$ \ and

\item {\it (globally) tame} if $\mathsf t (H) = \sup \{\mathsf t (H,u) \mid u \in \mathcal A (H_{\red}) \} < \infty$.
\end{itemize}

\smallskip
\noindent
{\bf Krull monoids.} We gather the basics of Krull monoids (detailed presentations can be found in \cite{HK98, Ge-HK06a}).
Let  $H$  and $D$  be  commutative and cancellative monoids and let $\varphi \colon H \to D$ be a monoid homomorphism. Then $\varphi$ can be uniquely lifted to a group homomorphism $\mathsf q ( \varphi ) \colon \mathsf q (H) \to \mathsf q (D)$,  where  $\mathsf q (H)$ and $\mathsf q (D)$ denote the quotient groups of $H$ and $D$. Moreover, $\varphi$  is called
\begin{itemize}
\smallskip
\item a  {\it divisor homomorphism} if $\varphi(a)\mid\varphi(b)$ implies that $a \t b$  for all $a,b \in H$;

\smallskip
\item  {\it cofinal} \ if for every $a \in D$ there exists some $u
      \in H$ such that $a \t \varphi(u)$;

\smallskip
\item  a {\it divisor theory} (for $H$) if $D = \mathcal F (P)$
for some set $P$, $\varphi$ is a divisor homomorphism, and for every
$p \in P$ (equivalently for every $a \in \mathcal{F}(P)$), there exists a finite
subset $\emptyset \ne X \subset H$ satisfying $p = \gcd \bigl(
\varphi(X) \bigr)$.
\end{itemize}
We call
$\mathcal{C}(\varphi)=\mathsf q (D)/ \mathsf q (\varphi(H))$ the
class group of $\varphi $ and use additive notation for this group.
For \ $a \in \mathsf q(D)$, we denote by \ $[a] = [a]_{\varphi} = a
\,\mathsf q(\varphi(H)) \in \mathsf q (D)/ \mathsf q (\varphi(H))$ \
the class containing \ $a$. The homomorphism $\varphi$ is cofinal if
and only if $\mathcal{C}(\varphi) = \{[a]\mid a \in D \}$, and if
$\varphi$ is a divisor homomorphism, then $\varphi(H)= \{a \in D
\mid [a]=[1]\}$. If $\varphi \colon H \to \mathcal F (P)$ is a
cofinal divisor homomorphism, then
\[
G_P = \{[p] = p \mathsf q (\varphi(H)) \mid p \in P \} \subset
\mathcal{C}(\varphi)
\]
is called the \ {\it  set of classes containing prime divisors}, and
we have $[G_P] = \mathcal{C}(\varphi)$.

Suppose that $H \subset D$ and that $\varphi = (H \hookrightarrow D)$.  Then, $\mathcal C ( \varphi ) = \mathsf q (D)/\mathsf q (H)$, and for $a \in D$ we set $[a]_{\varphi}=[a]_{D/H} = a \mathsf q (H)$. Then
\[
D/H = \{ [a]_{D/H} \mid a \in D\} \subset \mathcal C ( \varphi )
\]
is a submonoid with quotient group $\mathsf q (D/H)=\mathcal C (\varphi)$.  It is easy to check that $D/H$ is a group if and only if $H \hookrightarrow D$ is cofinal. In particular, if $D/H$ is finite or if $\mathsf q (D)/\mathsf q (H)$ is a torsion group, then $D/H= \mathsf q (D)/\mathsf q (H)$.

For subsets $A, B \subset \mathsf q (H)$, we denote by
$(A \DP B) = \{ x \in \mathsf q (H) \mid x B \subset A \}$, by $A^{-1} = (H \DP A)$, and by $A_v = (A^{-1})^{-1}$.  A subset $\mathfrak a \subset H$ is  an $s$-ideal of $H$ if $\mathfrak a H = \mathfrak a$. A subset $X \subset \mathsf q (H)$ is  a fractional  $v$-ideal (or a {\it fractional  divisorial ideal}) if there is a $c \in H$ such that $cX \subset H$ and $X_v = X$. We denote by $\mathcal F_v (H)$ the set of all fractional $v$-ideals and by $\mathcal I_v (H)$ the set of all $v$-ideals of $H$.
For $A, B \in \mathcal F_v (H)$, we define $A \cdot_v B = (A B)_v$, and we call $A \cdot_v B$ the $v$-product of $A$ and $B$.
Furthermore, $\mathcal I_v^* (H)$ is the monoid of $v$-invertible $v$-ideals (with $v$-multiplication) and $\mathcal F_v (H)^{\times} = \mathsf q \big( \mathcal I_v^* (H) \big)$ is its quotient group of fractional invertible $v$-ideals. The monoid $H$  is called $v$-noetherian if it satisfies the ascending chain condition on $v$-ideals. If $H$ is $v$-noetherian, then all sets of lengths are finite. We denote by $\mathfrak X (H)$ the set of all minimal non-empty prime $s$-ideals of $H$.

The map $\partial \colon H \to \mathcal I_v^* (H)$, defined by $\partial (a) = aH$ for each $a \in H$, is a cofinal divisor homomorphism. Thus, if  $\mathcal H = \{aH
\mid a \in H \}$ is the monoid of principal ideals of $H$, then we have a cofinal divisor homomorphism
\begin{equation} \label{v-class-group}
\mathcal H \hookrightarrow \mathcal I_v^* (H) \quad \text{and} \quad \mathcal C_v(H) = \mathcal I_v^* (H)/\mathcal H = \mathcal F_v (H)^{\times}/ \mathsf q (\mathcal H)
\end{equation}
is the {\it $v$-class group} of $H$. The monoid $H$ is called a {\it Krull monoid} if it satisfies one of
the following equivalent conditions (\cite[Theorem 2.4.8]{Ge-HK06a}){\rm \,:}
\begin{itemize}
\item[(a)] $H$ is $v$-noetherian and completely integrally closed.

\smallskip
\item[(b)] $H$ has a divisor theory.

\smallskip
\item[(c)] There is a free abelian monoid $F$ such that $H_{\red} \hookrightarrow F$ is a divisor homomorphism.
\end{itemize}
Let $H$ be a Krull monoid. Then a divisor theory
$\varphi \colon H \to \mathcal F (P)$ is unique up to unique
isomorphism. In particular, the class group $\mathcal C ( \varphi)$
defined via a divisor theory $\varphi$ of $H$ and the subset of classes
containing prime divisors depend only on $H$. Thus it is called the
{\it class group} of $H$ and is denoted by $\mathcal C (H)$. If $H$
is a Krull monoid, then $\mathcal{I}_v^{*}(H)$  is a free abelian monoid with basis
$\mathfrak X (H) = v$-$\spec (H) \setminus \{\emptyset \}$, the map $\partial \colon H\rightarrow
\mathcal{I}_v^{*}(H)$   is a divisor theory, and  $\mathcal{C}(H)$ is isomorphic to  the $v$-class group $\mathcal C_v (H)$ of $H$.

Main examples of Krull monoids stem from ring and module theory. We mention two examples and refer to \cite{Ge-HK06a, Ge-Zh20a} for more.   To begin with,   a commutative integral domain is  Krull  if and only its monoid of nonzero elements is a Krull monoid. Therefore, Condition (a) shows that every integrally closed noetherian domain is a Krull domain. To give an example from module theory, let $R$ be a ring and $\mathcal C$ be a small class of left $R$-modules that is closed under finite direct sums, direct summands, and isomorphisms. Then the set $\mathcal V (\mathcal C)$ of isomorphism classes of modules from $\mathcal C$ is a  commutative monoid, with operation induced by the direct sum. If the modules in $\mathcal C$ are not isomorphic to  proper direct summands, then $\mathcal V (\mathcal C)$ is unit-cancellative. If the endomorphism rings $\End_R (M)$ are semilocal for all modules $M$ from $\mathcal C$, then $\mathcal V (\mathcal C)$ is a Krull monoid (\cite[Theorem 3.4]{Fa02}, \cite{Fa19a}).

\smallskip
\noindent
{\bf Monoids of zero-sum sequences.} Let $G$ be an additive abelian group and let $G_0 \subset G$ be a subset. Then $[G_0] \subset G$ denotes the submonoid and $\langle G_0 \rangle \subset G$ denotes the subgroup generated by $G_0$. Let $\mathsf r (G)$ be the torsion-free rank of $G$.
In additive combinatorics, the elements of the free abelian monoid $\mathcal F (G_0)$ are called {\it sequences} over $G_0$ (they can be understood as finite unordered sequences of terms from $G_0$, with repetition allowed). For $S = g_1 \cdot \ldots \cdot g_{\ell} \in \mathcal F (G_0)$, let $\sigma (S) = g_1 + \ldots + g_{\ell} \in G$ denote the {\it sum} of $S$. Then the set
\[
\mathcal B (G_0 ) = \{ S \in \mathcal F (G_0) \mid \sigma (S) = 0 \} \subset \mathcal F (G_0)
\]
is a submonoid of $\mathcal F (G_0)$, called the {\it monoid of zero-sum sequences} over $G_0$. We say that $G_0$ is {\it condensed} if for every $g \in G_0$ there is $B \in \mathcal B (G_0)$ with $\mathsf v_g (B)>0$. Thus, $G_0$ is condensed if and only if the inclusion $\mathcal B (G_0) \hookrightarrow \mathcal F (G_0)$ is cofinal.  Since the inclusion $\mathcal B (G_0) \hookrightarrow \mathcal F (G_0)$ is a divisor homomorphism, $\mathcal B (G_0)$ is a Krull monoid by Condition (c) (given in the definition of Krull monoids). Moreover, by \cite[Proposition 2.5.6]{Ge-HK06a}, we have that
\begin{equation} \label{divisor-theory}
\mathcal B (G_0) \hookrightarrow \mathcal F (G_0) \ \text{is a divisor theory if and only if} \ \langle G_0 \rangle = [G_0 \setminus \{g\} ] \ \text{for every} \ g \in G_0 \,.
\end{equation}

If $G_0$ is condensed, then $\sigma\colon\mathcal{F}(G_0)\longrightarrow G$ induces by \cite[Proposition 2.5.62]{Ge-HK06a} an isomorphism
\[
\mathcal{C}(\block(G_0))\hookrightarrow \free(G_0))\longrightarrow [G_0].
\]
As usual, we set $\mathcal L (G_0) := \mathcal L \big( \mathcal B (G_0) \big)$, $\mathsf c (G_0) := \mathsf c \big( \mathcal B (G_0) \big)$, $\Delta (G_0) := \Delta \big( \mathcal B (G_0) \big)$, and so on.

\smallskip
\section{Transfer homomorphisms  of Krull monoids} \label{3}
\smallskip

We start with the definition of transfer homomorphisms in a not necessarily commutative setting (see \cite[Definition 2.1]{Ba-Sm15}). We need the concept in this general setting in Section \ref{4}.

\smallskip
\begin{definition} \label{3.1}
A monoid homomorphism $\theta \colon H \to B$  of atomic unit-cancellative monoids $H$ and $B$ is called a {\it transfer homomorphism} if it satisfies the following two properties:
\begin{enumerate}
\item[{\bf (T\,1)\,}] $B =  B^\times \theta(H) B^\times$ \ and \ $\theta ^{-1} (B^\times) = H^\times$.

\smallskip

\item[{\bf (T\,2)\,}] If $u \in H$, \ $b,\,c \in B$ \ and \ $\theta (u) = bc$, then there exist  $v,\,w \in H$ and $\varepsilon \in B^{\times}$ such that  $u = vw$,
                       $\theta (v)  =  b \varepsilon^{-1}$, and  $\theta (w)  = \varepsilon c $.
\end{enumerate}
\end{definition}

It is easy to check that the composition of transfer homomorphisms is a transfer homomorphism again. Moreover, transfer homomorphisms allow to pull back arithmetical properties from $B$ to $H$. In particular, if
$\theta \colon H \to B$ is a transfer homomorphism, then
\begin{equation} \label{transfer-lengths}
\mathsf L_H (a) = \mathsf L_B \bigl( \theta (a) \bigr) \quad \text{for every } \ a \in H \quad \text{and} \quad \mathcal L (H) = \mathcal L (B) \,.
\end{equation}
Details can be found in \cite[Section 2]{Ba-Sm15} and in \cite[Chapter 3.2]{Ge-HK06a}.
The next lemma describes the canonical transfer homomorphism from a Krull monoid to a monoid of zero-sum sequences (\cite[Chapter 3.4]{Ge-HK06a}).

\smallskip
\begin{lemma} \label{3.2}
Let $H$ be a Krull monoid, $\varphi \colon H \to F = \mathcal F (P)$ be a cofinal divisor homomorphism, $G = \mathcal C (\varphi)$ its class
group, and $G_P \subset G$ the set of classes containing prime
divisors. Let $\widetilde{\boldsymbol \beta} \colon F \to \mathcal F
(G_P)$ denote the unique homomorphism defined by
$\widetilde{\boldsymbol \beta} (p) = [p]$ for all $p \in P$.
\begin{enumerate}
\item The homomorphism $\boldsymbol \beta_H =  \boldsymbol \beta  = \widetilde{\boldsymbol \beta} \circ \varphi \colon H \to \mathcal B
      (G_P)$ is a transfer homomorphism.

\item The inclusion $\mathcal B (G_P) \hookrightarrow \mathcal F (G_P)$ is a cofinal divisor homomorphism.
\end{enumerate}
\end{lemma}

The homomorphism $\boldsymbol \beta$ in Lemma \ref{3.2} is called the {\it block homomorphism}, and $\mathcal B (G_P)$ is called the {\it block monoid}  associated  to $\varphi$. If $\varphi$ is a divisor theory, then $\mathcal B (G_P)$ is called the {\it block monoid} associated to $H$. Block homomorphisms enable us to study the arithmetic of Krull monoids with methods from additive combinatorics (see \cite[Chapter 1]{Ge-Ru09} and \cite{Sc16a, Fa-Tr18a} for surveys on this interaction).

If the inclusion $\mathcal B (G_P) \hookrightarrow \mathcal F (G_P)$ is a divisor theory, then $\varphi \colon H \to \mathcal F (P)$ is a divisor theory. However, in a variety of settings, the fact that the monoid under consideration is  Krull, is proved by establishing a (cofinal) divisor homomorphism to a free abelian monoid, which in general need not be a divisor theory. Moreover, oftentimes there is no information on the  class group of the Krull monoid but just on the class group of the divisor homomorphism (see \cite[Theorem 3.4]{Fa02} for an example from module theory and  \cite[Proposition 3.17]{Sm19a} for an example from ring theory).

The question, which pairs $(G, G_P)$ can actually occur, is answered by the following realization result (\cite[Theorems 2.5.4 and 3.7.8]{Ge-HK06a}). Let $G$ be an abelian group and $G_0 \subset G$ be a subset. Then the following statements are equivalent.
\begin{itemize}
\item[(a)] There is a Dedekind domain $D$ with class group $G$ such that $G_0$ is the set of classes containing non-zero prime ideals.

\item[(b)] There is a Krull monoid $H$ with class group $G$ such that $G_0$ is the set of classes containing prime divisors.

\item[(c)]  $G = [G_0]$.
\end{itemize}
Comparing this result with \eqref{divisor-theory} we see that, in general, the inclusion $\mathcal B (G_P) \hookrightarrow \mathcal F (G_P)$, given in Lemma \ref{3.2}, is not a divisor theory. By Condition (a), this need not even be the case when the homomorphism $\varphi$ (in Lemma \ref{3.2}) is a divisor theory. If $G_P$ consists of  torsion elements, then there is  the following result.

\smallskip
\begin{proposition} \label{3.3}
Let $H$ be a Krull monoid with torsion class group $G$. Then there is a subset $G_0^* \subset G$ such that $\mathcal B (G_0^*) \hookrightarrow \mathcal F (G_0^*)$ is a divisor theory and a transfer homomorphism $\theta^* \colon H \to \mathcal B (G_0^*)$.
\end{proposition}

\begin{proof}
Let $G_P \subset G$ denote the set of classes containing prime divisors.  By \cite[Corollary 4.13]{Sc10a}, there is a subset $G_0^* \subset G$ such that $g \in \langle G_0^* \setminus \{g\} \rangle$ for all $g \in G_0^*$ and a transfer homomorphism $\theta \colon \mathcal B (G_P) \to \mathcal B (G_0^*)$ (if $G_0$ is finite, this follows from a simple combinatorial argument given in \cite[Theorem 6.7.11]{Ge-HK06a}). Thus Equation \ref{divisor-theory} implies that $\mathcal B (G_0^*) \hookrightarrow \mathcal F (G_0^*)$ is a divisor theory. If $\boldsymbol \beta \colon H \to \mathcal B (G_P)$ is the transfer homomorphism given in Lemma \ref{3.2}, then $\theta^* = \theta \circ \boldsymbol \beta \colon H \to \mathcal B (G_0^*)$ is a transfer homomorphism too.
\end{proof}

The goal of this section is to obtain a result as given in  Proposition \ref{3.3} for Krull monoids with finitely generated class group. Starting from cases where the torsion-free rank of the class group equals one (\cite{B-C-R-S-S10, Ba-Kr11a, Ge-Gr-Sc-Sc10, De-Ze20a}), the arithmetic of Krull monoids with finitely generated class group has seen strong interest (e.g., \cite{Ga-La-Sm19a, Gr22a}), partly motivated by module theory (\cite{Fa06a, Fa12a, Ba-Ge14b, Fa19a,Ca21a}). Moreover, the refinement of the block homomorphism, as given in Proposition \ref{3.3} and Theorem \ref{T:MainThm}, is also used for studying the arithmetic of Krull monoids which have prime divisors in all classes (see Remark \ref{3.9}). We formulate our main result.

\smallskip
\begin{theorem} \label{T:MainThm}
Let $H$ be a Krull monoid with finitely generated class group. Then there exist a finitely
generated abelian group $G^*$,  a subset $G_0^* \subset G^*$ such that $\block(G_0^*)\hookrightarrow \free(G_0^*)$ is a
divisor theory, and a transfer homomorphism
$\theta^* \colon H \rightarrow \block(G_0^*)$, which is  a composition of block homomorphisms. In particular, we have
\begin{enumerate}
\item $\mathcal L (H) = \mathcal L (G_0^*)$.

\item $\mathsf c (G_0^*) \le \mathsf c (H) \le \max \{ \mathsf c (G_0^*), 2 \}$.

\item $H$ is locally tame if and only if $\mathcal B (G_0^*)$ is locally tame.
\end{enumerate}
\end{theorem}

\smallskip
To prove Theorem \ref{T:MainThm}, we need some preparations.
First  we construct, for a given reduced Krull monoid $H$ and a divisor homomorphism to a free abelian monoid, a divisor theory for $H$. Such a construction has already been given in \cite[Theorem 2.4.7]{Ge-HK06a}, but this
construction is not explicit enough for our purpose. So, let $H$ be a reduced Krull monoid and  $F=\free(P)$ be a free abelian monoid such that the inclusion
$H\hookrightarrow F$ is a divisor homomorphism.
We say that an element  $x\in F$ is irrelevant if $(xF\cap H)_v=H$.

\smallskip
\begin{lemma}\label{L:vPrimes}\hfill
\begin{enumerate}
\item The product of irrelevant elements is irrelevant, too.

\item Let $p\in P$ be not irrelevant. Then $pF\cap H\in \vspec(H)$.
\item Let $q\in \mathfrak X (H)$. Then there exists some $p\in P$, such that $q=pF\cap
H$.
\end{enumerate}
\end{lemma}

\begin{proof}
1.  If $x$, $y\in F$ are irrelevant, then
\[
H=H\cdv H=(xF\cap H)_v\cdv (yF\cap H)_v=\Big((xF\cap H)(yF\cap H)\Big)_v\subset (xyF\cap
H)_v\subset H \,.
\]

2. If $pF\cap H=\emptyset$, we are done. Hence we may assume $pF\cap H\not=\emptyset$. Then $pF\cap
H$ is a non-empty prime $s$-ideal of $H$, and hence contains some minimal, non-empty, prime
$s$-ideal $q$. Then $q\in\mathfrak X (H)$. Now we have
\[
q\subset pF\cap H\subset (pF\cap H)_v\not=H  \,.
\]
Since $H$ is Krull, $q$ is a maximal $v$-ideal, which implies $q=pF\cap H$.

\smallskip
3. Let $x=\gcd(q)\in F$. By \cite[Proposition 2.4.2.6]{Ge-HK06a} we have $q=xF\cap H$. Therefore
$x$ is not irrelevant. Now let $p_1$,\ldots , $p_n\in P$ be such that $x=p_1\ldots p_n$. By Lemma
\ref{L:vPrimes}.1 we may assume that $p_1$ is not irrelevant. Then have
\[
q=xF\cap H\subset p_1F\cap H  \,.
\]
By 1. we know, that $p_1F\cap H$ is a $v$-ideal, not equal to $H$. Since  $q$ is $v$-maximal, we
obtain $q=p_1F\cap H$.
\end{proof}

\medskip
By Lemma \ref{L:vPrimes}.3, we may choose a map $f\colon \mathfrak X (H)\rightarrow P$ such that $q=f(q)F\cap H$ for each $q\in\mathfrak X (H)$ (in particular $f$ is injective).
Now let $q\in \mathfrak X (H)$ and let
\[
\mathsf v_q\colon \mathsf q(H)\rightarrow \Z,\quad \mathsf v_{f(q)}\colon \mathsf q(F)\rightarrow \Z
\]
be the corresponding valuations. From $q=f(q)F\cap H$ we obtain $H_q^\times =H_q\cap
F_{f(q)F}^\times$. Hence, if $\pi$ is a prime of $H_q$ (which is a discrete valuation monoid), we have $\pi=\ve
f(q)^{e(q)}$ for some $\ve\in F_{f(q)F}^\times$ and some $e(q)\in\N$. This implies
\[
\mathsf v_{f(q)}(h) = e(q) \mathsf v_q(h)
\]
for all $h\in \mathsf q(H)$. Therefore, we obtain that
\[
\mathsf v_{f(q)}(H) = e(q) \mathsf v_q(H)=e(q)\N_0
\]
and
\begin{equation}\label{E:eq}
e(q) = \min \mathsf v_{f(q)}(H)\setminus \{ 0\} = \gcd( \mathsf v_{f(q)}(H)) \,.
\end{equation}
We define
\begin{gather*}
P_0=f(\mathfrak X (H))\subset P,\quad F_0=[\{ f(q)^{e(q)}\st q\in
\mathfrak X (H)\}\cup P\setminus P_0], \quad \text{and} \\
F_1=[\{ f(q)^{e(q)}\st q\in
\mathfrak X (H)\}]  \,.
\end{gather*}
By \eqref{E:eq}, we have $H\subset F_0$. Since $H \hookrightarrow F$ is a divisor homomorphism, the same is true for  $H \hookrightarrow F_0$.
We define $d \colon F_0\rightarrow F_1$ by
\[
d_{\mid F_1}=\id_{F_1},\quad  d(p)=1,\, \text{ for all } p\in P\setminus P_0
\]
and set $\partial=d_{\mid H}\colon H\rightarrow F_1$.

Note that $F_1$ is free with basis $\{f(q)^{e(q)}\st q\in\mathfrak X (H)\}$. Hence  $f$ induces an
isomorphism $f_* \colon \DivMon(H)\rightarrow F_1$, such that $f_*(q)=f(q)^{e(q)}$ for all
$q\in\mathfrak X (H)$.

By construction of the  $e(q)$'s the following diagram is commutative
\[
\begin{xy}
\qtriangle/{->}`{->}`{<-}/[{H}`{F_1}`{\DivMon(H)};{\partial}`{\partial_H}`{f_*}]
\end{xy}   \,,
\]
where $\partial_H$ is the canonical divisor theory of $H$, defined by $\partial_H(a)=aH$ for all
$a\in H$. Therefore, $\partial$ is a divisor theory of $H$, and $f_*$ induces an isomorphism
\[
\mathcal C(\partial)=\coker( \mathsf q(\partial))\cong \coker( \mathsf q(\partial_H))= \mathcal C_v(H)  \,.
\]

Let $G = \mathsf q(F) / \mathsf q(H)$ be the class group of $H\subset F$, and set $\overline{G}= \mathsf q(F_0)/ \mathsf q(H)\subset G$.
Obviously, we have
\[
\mathsf q(F)/ \mathsf q(F_0)\cong \bigoplus_{q\in\mathfrak X (H)}\Z/e(q)\Z   \,.
\]
Hence, we obtain an exact sequence
\begin{equation}\label{E:ExSeq}
0\rightarrow \overline{G}\rightarrow
G\rightarrow \bigoplus_{q\in\mathfrak X (H)}\Z/e(q)\Z\rightarrow 0 \,.
\end{equation}
Since $\partial=d_{\mid H}$, $\mathsf q (d)$ induces an epimorphism
\[
\delta\colon \overline{G}= \mathsf q(F_0)/ \mathsf q(H)\rightarrow \coker(\mathsf q(\partial))\cong \mathcal C_v(H) \,.
\]
Thus we obtain the following result.

\medskip
\begin{proposition} \label{P:EqualRank}
Let all notation be as above.
\begin{enumerate}
\item $\mathcal C_v(H)$ is isomorphic to a subquotient of the class group $G$ of $H \hookrightarrow F$. In particular,  we have
      \[
      \mathsf r(G) \ge \mathsf  r( \mathcal C_v(H)) \, .
      \]

\item If $\mathsf r(G) = \mathsf r(\mathcal C_v(H))<\infty$, then $P_0=P$, $F_1=F_0$, and $\overline{G}\cong \mathcal C_v(H)$.
\end{enumerate}
\end{proposition}

\begin{proof}
1. This follows immediately from the above construction.

2.   From the exact sequence \eqref{E:ExSeq} we get $\mathsf{r}(\overline{G})=\mathsf{r}(G)$. Hence $\mathsf{r}(\overline{G})=\mathsf{r}(\mathcal C_v(H))$ by assumption.

We consider now the epimorphism $\delta\colon \overline{G}\rightarrow \mathcal C_v(H)$ from above.
Since $\mathsf r(\mathcal C_v(H)) =
\mathsf r(\overline{G})<\infty$,  $\delta\otimes\Q$ is an
isomorphism, which implies that the kernel of $\delta$ consists  only of elements with finite order.

Assume to the contrary that $P_0\not=P$, and let $p\in P\setminus P_0$. If $g$ is the image of $p$ in $\overline{G}$,
we have by definition $\delta(g)=0$. Hence $g$ is of finite order, which means $p^n\in H$ for some
$n\in\N$. From $H^\times =H\cap F^\times$, we obtain $p^nH\not=H$. Since $H \hookrightarrow F$ is
a divisor homomorphism, we have $p^nH=p^nF\cap H$. Therefore $p^n$ is not irrelevant. From Lemma \ref{L:vPrimes}.1 we
deduce that $p$ is not irrelevant. Since $p^n\in pF\cap H$, Lemma \ref{L:vPrimes}.2 tells us
$q:=pF\cap H\in \vspec(H)\setminus\{ \emptyset\}=\mathfrak X (H)$. From
\[
p^n\in q=f(q)F\cap H
\]
we obtain $p=f(q)\in P_0$, a contradiction. Now $P_0=P$ implies $F_0=F_1$, $\partial= H\hookrightarrow F_0$ and $\bar{G}=\coker(q(\partial))\cong \mathcal C_v(H)$.
\end{proof}

We apply  this result now for block monoids.
Let $G$ be an abelian group, $G_0\subset G$ a condensed subset.  We assume that $\mathsf{r}(\mathcal{C}_v(\block(G_0)))=\mathsf{r}([G_0])$. For $g\in G_0$ we set
\[
e_{G_0}(g)=\min \mathsf{v}_g(\block(G_0))\setminus\{ 0\}=\gcd\mathsf{v}_g(\block(G_0)).
\]
and define $G_0^1=\{e_{G_0}(g)g\st g\in G_0\}$. Then in the notation from \ref{P:EqualRank} we have $\block(G_0)\subset F_1=[\{g^{e_{G_0}(g)}\}]$ and this inclusion $\iota$  is a divisor theory. Since $G_0$ is condensed, the sum homomorphism $\sigma\colon\free(G_0)\longrightarrow G$ induces an isomorphism $\free(G_0)/\block(G_0)\longrightarrow [G_0]$. By restriction we get an isomorphism $\bar{\sigma}\colon F_1/\block(G_0)\longrightarrow [G_0^1]$ with
\[
\bar{\sigma}([\{g^{e_{G_0}(g)}]\st g\in G_0\})=G_0^1.
\]
By identifying the class group of $\iota$ with $[G_0^1]$ by means of $\bar{\sigma}$ we obtain the block homomorphism
\[
\beta_{G_0}\colon\block(G_0)\longrightarrow\block(G_0^{1})
\]
defined by
\[
\prod_{g\in G_0} g^{a_g}=\prod_{g\in G_0} (g^{e_{G_0}(g)})^{a_g/e_{G_0}(g)}\mapsto \prod_{g\in G_0} (e_{G_0}(g)g)^{a_g/e_{G_0}(g)}.
\]
\begin{lemma}\label{L:DivTh}
In the above situation the following are equivalent:
\begin{enumerate}
\item $\block(G_0)\hookrightarrow \free(G_0)$ is a divisor theory.
\item For all $g\in G_0$ we have $e_{G_0}(g)=1$.
\item $[G_0]=[G_0^1]$.
\end{enumerate}
\end{lemma}
\begin{proof}
1 $\Longrightarrow$ 2: so suppose  that $\block(G_0)\hookrightarrow \free(G_0)$ is a divisor theory and let $g\in G_0$. Then $g=\gcd(X)$ for some $\emptyset\not=X\subset\block(G_0)$. In particular we have $v_g(B)=1$ for some $B\in \block(G_0)$. Hence $e_{G_0}(g)=1$.

2 $\Longrightarrow$ 3: by our identifications of the class groups of the inclusions $\block(G_0)\hookrightarrow F_1$, $\block(G_0)\hookrightarrow \free(G_0)$ the exact sequence \eqref{E:ExSeq} translates into the exact sequence
\[
0\longrightarrow [G_0^1]\longrightarrow [G_0] \longrightarrow \bigoplus_{g\in G_0}\Z/e_{G_0}(g)\Z\longrightarrow 0\, .
\]

3 $\Rightarrow$ 1: If $[G_0^1]=[G_0]$ we get by the identifications above $q(F_1)/q(\block(G_0))=q(\free(G_0))/q(\block(G_0))$. Therefore $q(F_1)=q(F_0)$. Looking at the the definition of $F_1$ one immediately deduces $\free(G_0)=F_1$. Since $\block(G_0)\hookrightarrow F_1$ is a divisor theory, we are ready.
\end{proof}

\medskip
\begin{proof}[Proof of Theorem \ref{T:MainThm}]
Let $H$ be reduced Krull monoid with finitely generated class group.
Choose a divisor theory $H\hookrightarrow F = \mathcal F (P)$   with class group $G$. Let
$\pi\colon F\rightarrow G$ be the canonical epimorphism, and set $G_0=\pi(P)\subset G$. By Lemma \ref{3.2},
$G_0\subset G$ is a condensed subset,  $\mathcal B (G_0) \hookrightarrow \mathcal F (G_0)$ is a cofinal divisor homomorphism, and we consider the block homomorphism $\boldsymbol \beta_H \colon H \to \mathcal B (H) := \mathcal B (G_0)$.
Since $\block(G_0)$ is again reduced and  Krull, we can iterate this construction. Indeed,  for $n\in\N_0$ we define inductively
\begin{equation} \label{iteration}
\block^0(H)=H,\quad \block^{n+1}(H)=\block(\block^n(H)),\quad \boldsymbol \beta^0_{H}=\id_{H}\colon
H\rightarrow\block^0(H),\quad \boldsymbol \beta^{n+1}_H=\boldsymbol \beta_{\block^n_H}\circ \boldsymbol \beta^n_H\colon H\rightarrow
\block^{n+1}_H    \,.
\end{equation}
For $n\in\N$, let $\overline{G}^n$ be the class group of $\block^{n-1}(H)$. Then, by construction, we have
\[
\block^n(H)=\block(\overline{G}^n_0)
\]
for a condensed subset $\overline{G}^n_0\subset \overline{G}^n$ such that $[\overline{G}^n_0]=\overline{G}^n$. It follows from Proposition
\ref{P:EqualRank}.1 that we have a descending sequence of ranks:
\begin{equation}\label{E:Ranks}
\mathsf r(\overline{G}^1) \geq \mathsf r(\overline{G}^2) \geq \mathsf r(\overline{G}^3)\geq\ldots \,.
\end{equation}
and that all $\overline{G}^n$ are finitely generated, too. The  sequence of ranks stabilizes. We choose $k\in\N$ such $\mathsf{r}(\overline{G}^{n+1})=\mathsf{r}(\overline{G}^n)$ for all $n\geq k$ and set $G^n=\overline{G}^{n+k}$, $G^n_0=\overline{G}_0^{k+n}$ for all $n\in\N_0$. Then we can assume that for all $n\in\N_0$  the block homomorphism  $\block(G_0^n)\longrightarrow\block(G_0^{n+1})$ is the block homomorphism $\beta_{G_0^n}$ constructed above. In particular we have
\[
G^0=[G_0^0]\supset G^1=[G_0^1]\supset\ldots\quad\text{and}\quad  G_0^{n+1}=\{e_{G_0^n}(g)g\st g\in G_0^n\}
\]
for all $n\in\N_0$.

For all $n\in\N_0$ we define
\[
\gamma^{n,n+1}\colon G_0^n\longrightarrow G_0^{n+1},\quad g\mapsto e_{G_0^n}(g)g,\quad \gamma^n=\gamma^{n,n-1}\circ\ldots\circ\gamma^{0,1}\colon G_0^0\longrightarrow G_0^n
\]
(in particular $\gamma^0=\id_{G^0_0}$). For $n\in\N_0$ and $g\in G_0^0$ we define inductively
\[
e^0(g)=1,\quad e^{n+1}(g)=e_{G_0^n}(\gamma^n(g))e^n(g).
\]
Then an easy induction shows for all $n\in\N_0$:
\begin{itemize}
\item for  all $g\in G_0$ we have $\gamma^n(g)=e^n(g)g$;
\item the composition of block homomorphisms
\[
\beta^n\colon\block(G_0^0)\overset{\beta_{G_0^0}}{\longrightarrow}\block(G_0^1)\overset{\beta_{G_0^0}}{\longrightarrow}\ldots \overset{\beta_{G_0^{n-1}}}{\longrightarrow}\block(G_0^n)
\]
is given as follows: if $B=\prod_{g\in G_0}g^{a_g}\in\block(G_0)$ then
\[
\beta^n(B)=\prod_{h\in G_0^n}h^{b_h}\, ,
\]
where for each $h\in G_0^n$
\[
b_h=\sum_{\substack{g\in G_0\\\gamma^n(g)=h}}\frac{a_g}{e^n(g)}.
\]
In particular
\[
\vert\beta^n(B)\vert=\sum_{g\in G^0_0}\frac{a_g}{e^n(g)}.
\]
\end{itemize}

Now let $B\in\block(G_0^0)$. Since for each $g\in G_0^0$ the sequence $(e^n(g))_{n\in\N_0}$ is increasing, we get the decreasing sequence of positve integers
\[
\vert B\vert\geq \vert\beta^1(B)\vert\geq \ldots .
\]
This sequence has to stabilize. Using again, that $(e^n(g))_{n\in\N_0}$ is increasing, we see that for each $g\in\supp(B)=\{ g\in G_0^0\st a_g>0\}$, that the sequence $(e^n(g))_{n\in\N_0}$ is bounded. Since any $g\in G_0^0$ is contained in some $\supp(B)$, $B\in\block(G_0)$, we see that $(e^n(g))_{n\in\N_0}$ is bounded for any $g\in G_0^0$. We set $e^{\infty}(g)=\max\{e^n(g)\st n\in\N_0\}$. Since $e^0(g)\mid e^1(g)\ldots $ we get $e^n(g)\mid e^{\infty}(g)$ for any $n\in\N$ and any $g\in G_0^n$.

We define $G^\infty=\langle\{ e^{\infty}_{G_0}g\st g\in G_0\}\rangle$. Since $G^0=[G_0]$, the group
$G^0/G^{\infty}$ is a finitely generated torsion group, hence finite.

We obtain in $G^0$ the descending sequence of subgroups
\[
G^0\supset G^1\supset G^2\supset\ldots \supset G^{\infty}  \,.
\]
Since $G^0/G^{\infty}$ is finite, we have $[G^n_0]=G^n=G^{n+1}=[G^{n+1}_0]$ for some $n\in\N$. Using \ref{L:DivTh} we see that $\block(G^n_0)\hookrightarrow \free(G^n_0)$ is a divisor theory.

\smallskip

It remains to verify the statements 1. -- 3.

1. This holds true by \eqref{transfer-lengths}.

2. The asserted inequalities hold for  block homomorphisms by \cite[Theorem 3.4.10.2]{Ge-HK06a}, whence they hold for compositions of block homomorphisms.

3. For a Krull monoid $S$ with block homomorphism $\boldsymbol \beta_S$, $S$ is locally tame if and only if $\boldsymbol \beta_S (S)$ is locally tame by \cite[Proposition 3.3]{Ga-Ge-Sc15a}. Since $\theta^*$ is a composition of block homomorphisms, it follows that $H$ is locally tame if and only if $\mathcal B (G_0^*)$ is locally tame.
\end{proof}

\smallskip
\begin{remark}[{\bf Uniqueness of transfer homomorphisms}] \label{3.7}
Let $H$ be a Krull monoid. The divisor theory of $H$ is uniquely determined. Suppose that the class group of $H$ is torsion or finitely generated. Then Proposition \ref{3.3} and Theorem \ref{T:MainThm} imply that there is a subset $G_0^*$ of an abelian group and a transfer homomorphism $\theta^* \colon H \to \mathcal B (G_0^*)$ such that the inclusion $\mathcal B (G_0^*) \hookrightarrow \mathcal F (G_0^*)$ is a divisor theory. In contrast to the divisor theory of $H$, $G_0^*$ and $\theta^*$ are not uniquely determined in general. To give a simple example, consider a Krull monoid $H$ with infinite cyclic class group $G$, say $G = \langle e \rangle$, and suppose that $G_P = \{-2e, -e, 0,  e, 2e \} \subset G$ is the set of classes containing prime divisors. By Lemma \ref{3.2}, $\boldsymbol \beta_H \colon H \to \mathcal B (G_P)$ is a transfer homomorphism, and by \eqref{divisor-theory} the inclusion $\mathcal B (G_P) \hookrightarrow \mathcal F (G_P)$ is a divisor theory. By \cite[Proposition 6.12]{Ba-Ge14b}, there is a transfer homomorphism $\theta \colon \mathcal B (G_P) \to \mathcal B (C_3)$, where $C_3$ is a cyclic group with three elements. Thus, $\theta^* = \theta \circ \boldsymbol \beta_H \colon H \to \mathcal B (C_3)$ is a transfer homomorphism and the inclusion $\mathcal B (C_3) \hookrightarrow \mathcal F (C_3)$ is a divisor theory.

This example shows in particular that the iteration process, given in \eqref{iteration}, need not stabilize when the inclusion $\mathcal B (G_0^n) \hookrightarrow \mathcal F (G_0^n)$ is a divisor theory for the first time. The question, when this process stabilizes, was investigated in detail in \cite{Sc10a}, in case of torsion class groups.
\end{remark}

\smallskip
\begin{remark}[{\bf On global tameness}] \label{3.8}
Let $H$ be a Krull monoid and let $\boldsymbol \beta_H \colon H \to \mathcal B (G_P)$ be the transfer homomorphism, given in Lemma \ref{3.2}. Then $H$ is locally tame if and only if $\mathcal B (G_P)$ is locally tame. If the class group is finitely generated, then $H$ is tame if and only if $G_P$ is finite; but this equivalence does not hold in general (\cite[Theorem 4.2 and Example 4.13]{Ge-Ka10a}).
\end{remark}

\smallskip
\begin{remark}[{\bf On the set of minimal distances}] \label{3.9}
Let $H$ be an atomic monoid. The {\it set of minimal distances} $\Delta^* (H)$ is defined as the set of all $\min \Delta (S)$ over all divisor-closed submonoid $S \subset H$ with $\Delta (S) \ne \emptyset$. It plays a crucial role in structural descriptions of sets of lengths (see \cite[Chapter 4]{Ge-HK06a}). Suppose that $H$ is a Krull monoid with class group $G$ and suppose that each class contains a prime divisor. By \cite[Proposition 4.3.13]{Ge-HK06a}, we have
\[
\Delta^* (H) = \{ \min \Delta (G_0) \mid G_0 \subset G, \ \Delta (G_0) \ne \emptyset \} \,.
\]
Thus, when studying the set of minimal distances $\Delta^* (H)$ of a Krull monoid with class group $G$ and prime divisors in all classes, we need to study the set of distances $\Delta (G_0)$ of Krull monoids $\mathcal  B (G_0)$, and for doing so we use Theorem \ref{T:MainThm} (see \cite[Chapter 6.7]{Ge-HK06a}, \cite{Pl-Sc20a, Ge-Zh16a}).
\end{remark}

\smallskip
\section{Transfer Krull monoids and weakly Krull monoids} \label{4}
\smallskip

Transfer Krull monoids are monoids that allow a (suitably defined) transfer homomorphism onto a monoid of zero-sum sequences. By Lemma \ref{3.2}, Krull monoids are transfer Krull, but transfer Krull monoids need neither be commutative, nor cancellative, nor completely integrally closed, nor $v$-noetherian. We start with the definition of a weak transfer homomorphism as given in \cite[Definition 2.1]{Ba-Sm15}.

\smallskip
\begin{definition} \label{4.1}
Let $H$ be a monoid.
\begin{enumerate}
\item A monoid homomorphism $\theta \colon H \to B$ between  atomic unit-cancellative monoids $H$ and $B$ is called a {\it weak transfer homomorphism} if it has the following two properties.
    \begin{itemize}
    \item[{\bf (T1)}] $B = B^{\times} \theta (H) B^{\times}$ and $\theta^{-1} (B^{\times})=H^{\times}$.

    \item[{\bf (WT2)}] If $a \in H$, $n \in \N$, $v_1, \ldots, v_n \in \mathcal A (B)$ and $\theta (a) = v_1 \cdot \ldots \cdot v_n$, then there exist $u_1, \ldots, u_n \in \mathcal A (H)$ and a permutation $\tau \in \mathfrak S_n$ such that $a = u_1 \cdot \ldots \cdot u_n$ and $\theta (u_i) \in B^{\times} v_{\tau (i)} B^{\times}$ for each $i \in [1,n]$.
    \end{itemize}

\smallskip
\item $H$ is said to be a {\it transfer Krull monoid} (over $G_0)$ if $H$ is atomic and unit-cancellative and if there exists a weak transfer homomorphism $\theta \colon H \to \mathcal B (G_0)$ for a subset $G_0$ of an abelian group $G$. A domain is said to be transfer Krull if its monoid of regular elements is a transfer Krull monoid.
\end{enumerate}
\end{definition}

Let $\theta \colon H \to B$ be a homomorphism of atomic unit-cancellative monoids. If  $\theta $ is a transfer homomorphism, then $\theta$ is a weak transfer homomorphism and the converse holds if $H$ and $B$ are both commutative (\cite[Section 2]{Ba-Sm15}). However, there are monoids which do not have a transfer homomorphism to any commutative monoid but which do have a weak transfer homomorphism to a Krull monoid (\cite[Remark 2.4]{Ba-Sm15} and \cite[Example 4.5]{Ba-Ba-Go14}).

\smallskip
\begin{example} \label{4.2}~

\noindent
We mention some key examples of transfer Krull monoids that are not Krull, and refer to  \cite[Section 5]{Ge-Zh20a} for more.

1. Let $R$ be a bounded hereditary noetherian prime ring. If all stably free right $R$-ideals are free, then $R$ is a transfer Krull domain (\cite[Theorem 4.4]{Sm19a}). For more of this flavor see \cite{Ba-Sm15}.

2. A Mori domain is Krull if and only if it completely integrally closed. There are results stating that Mori domains, that are close to their complete integral closure, are transfer Krull (for a sample, see \cite[Theorem 5.8]{Ge-Ka-Re15a}).

3. Let $R$ be a Bass ring and $\mathcal C$ be the class of finitely generated torsion-free $R$-modules. Then $\mathcal V (\mathcal C)$ is a commutative, but not necessarily cancellative  transfer Krull monoid over a subset of a finitely generated abelian group (\cite[Theorem 1.1]{Ba-Sm21a}).
\end{example}

\smallskip
\begin{theorem} \label{4.3}
Let $H$ be a transfer Krull monoid over a subset of an abelian group $G$ that is  finitely generated resp. torsion. Then there exist an abelian group $G^*$ that is  finitely generated resp. torsion,  a subset $G_0^* \subset G^*$ such that $\block(G_0^*)\hookrightarrow \free(G_0^*)$ is a
divisor theory, and a weak transfer homomorphism
$\theta^* \colon H \rightarrow \block(G_0^*)$. In particular, we have
\begin{enumerate}
\item $\mathcal L (H) = \mathcal L (G_0^*)$.

\item If $H$ is commutative, then $H$ has finite catenary degree if and only if $\mathcal B (G_0^*)$ has finite catenary degree.
\end{enumerate}
\end{theorem}

\noindent
{\it Remark.} For cancellative but not necessarily commutative monoids $H$, axiomatically defined distance functions $\mathsf d$ and associated catenary degrees $\mathsf c_{\mathsf d} (H)$ are studied in \cite[Section 4]{Ba-Sm15} and in \cite{Sm19a}. The finiteness of these catenary degrees is also preserved under transfer homomorphisms, as it is the case in the commutative setting.

\begin{proof}
Let $G$ be an abelian group that is  finitely generated resp. torsion, $G_0 \subset G$ be a subset, and $\theta_1 \colon H \to \mathcal B (G_0)$ be a weak transfer homomorphism. By Theorem \ref{T:MainThm}, resp. by Proposition \ref{3.3},  there exist a finitely
generated abelian group $G^*$, resp. a torsion group $G^*$,  a subset $G_0^* \subset G^*$ such that $\block(G_0^*)\hookrightarrow \free(G_0^*)$ is a
divisor theory, and a transfer homomorphism
$\theta_2 \colon \mathcal B (G_0) \rightarrow \block(G_0^*)$. Then $\theta^* = \theta_2 \circ \theta_1 \colon H \to \mathcal B (G_0^*)$ is a weak transfer homomorphism.
Thus, it remains to verify the in particular statements.

As in the case of transfer homomorphisms, $\mathcal L (H) = \mathcal L (G_0^*)$ follows easily from the definition. Now suppose that $H$ is commutative. Since $H$ is unit-cancellative, the proof runs along the same lines as in the cancellative setting (\cite[Theorem 3.2.5]{Ge-HK06a}).
\end{proof}

\smallskip
Next we consider weakly Krull monoids and weakly Krull domains. Weakly Krull domains were introduced by Anderson, Anderson, Mott, and Zafrullah (\cite{An-An-Za92b, An-Mo-Za92}) and Halter-Koch showed their purely multiplicative character (\cite{HK95a}). We refer to \cite{HK98} for the ideal theory of weakly Krull monoids, to \cite[Example 5.7]{Ge-Ka-Re15a} for an extended list of examples, and to \cite{Fa-Wi21a} for weakly Krull monoid algebras. We recall the definition.  A commutative and cancellative monoid is
\begin{itemize}
\item  {\it weakly Krull} (\cite[Corollary 22.5]{HK98}) if
       \[
       H = \bigcap_{\mathfrak p \in \mathfrak X (H)} H_{\mathfrak p} \quad \text{and} \quad \{ \mathfrak p \in \mathfrak X (H) \mid a \in \mathfrak p \} \ \text{is finite for all} \ a \in H \,.
       \]

\item {\it weakly factorial} if one of the following equivalent conditions is satisfied (\cite[Exercise 22.5]{HK98}):
       \begin{itemize}
       \item Every non-unit is a finite product of primary elements.
       \item $H$ is a weakly Krull monoid with trivial $t$-class group.
       \end{itemize}
\end{itemize}
A commutative  domain $D$ is weakly Krull resp. weakly factorial if its monoid of nonzero elements if weakly Krull resp. weakly factorial.
The arithmetic of weakly Krull monoids is studied via $T$-block monoids (\cite[Chapter 3.4]{Ge-HK06a}). To recall this concept, let $G$ be an additively abelian group, $G_0 \subset G$ a subset, $T$ a commutative and cancellative monoid, $\sigma \colon \mathcal F (G_0) \to G$, the sum, and $\iota \colon G \to G$ a monoid homomorphism. Then
\[
B = \mathcal B (G_0, T, \iota) = \{ St \in \mathcal F (G_0)\times T \mid \sigma (S) + \iota (t) = 0 \} \subset \mathcal F (G_0) \times T = F
\]
is the {\it $T$-block monoid over $G_0$ defined by $\iota$}. The inclusion $B \hookrightarrow F$ is a divisor homomorphism, and if $G_0=G$ or if $G_0 \cup \iota (T)$ consists of torsion elements, then $B \hookrightarrow F$ is cofinal.

\medskip
\begin{proposition} \label{4.4}
Let $D = \mathcal F (P) \times T$ be a reduced atomic monoid, where $P \subset D$ a set of primes
and $T \subset D$ is a  submonoid. Let $H \subset D$ be a submonoid such that $H \hookrightarrow D$ is a cofinal divisor homomorphism with
class group $G = \mathsf q (D)/\mathsf q (H)$, and let $ G_P = \{[p]= p \mathsf q (H) \mid p \in P\} \subset G$ be the set of classes containing primes.
Let $\iota \colon T \to G$ be defined by $\iota (t) = [t]$, \ $F= \mathcal F (G_P) \times T$, \ $B =
\mathcal B(G_P,T, \iota) \subset F$, and let $\widetilde{\boldsymbol \beta} \colon D \to F$ be the
unique homomorphism satisfying $\widetilde{\boldsymbol \beta} (p) = [p]$ for all $p \in P$ and
$\widetilde{\boldsymbol \beta} \t T = \text{\rm id}_T$.

\begin{enumerate}
\item The restriction $\boldsymbol \beta = \widetilde{\boldsymbol \beta} \t H \colon H \to B$ \ is a transfer homomorphism satisfying $\mathsf c (H, \boldsymbol \beta) \le 2$.

\item  $B \hookrightarrow F$ is a cofinal divisor homomorphism, and there is an isomorphism $\overline \psi \colon \mathsf q (F)/\mathsf q (B) \to G$,
such that $\overline{\psi}(S\, t)=\sigma(S)+\iota(t)$ for all $S\, t\in \mathcal F(G_P) \times T$.
\end{enumerate}
\end{proposition}

\begin{proof}
1. and 2. See \cite[Proposition 3.4.8]{Ge-HK06a}.
\end{proof}

\smallskip
\begin{theorem} \label{4.5}
Let $H$ be a $v$-noetherian weakly Krull monoid with $\emptyset \ne \mathfrak f = (H \colon \widehat H) \subsetneq H$ and suppose that $\mathcal C_v (H)$ is finite. Let
\[
\mathcal P^* = \{\mathfrak p \in \mathfrak X (H) \mid \mathfrak p \supset \mathfrak f \}, \ \mathcal P = \mathfrak X (H) \setminus \mathcal P^* , \quad \text{and} \quad T = \prod_{\mathfrak p \in \mathcal P^*} (H_{\mathfrak p})_{\red} \,.
\]
\begin{enumerate}
\item There is a transfer homomorphism  $\boldsymbol \beta \colon H \to \mathcal B (G_{\mathcal P}, T, \iota)$, where   $G_{\mathcal P} \subset \mathcal C_v (H)$ is the set of classes containing minimal prime ideals. Moreover, $\mathcal B (G_{\mathcal P}, T, \iota)$ is a $v$-noetherian weakly Krull monoid.

\item For every divisor-closed submonoid $S \subset H$ with $\mathsf v_{\mathfrak p} (aH)=0$ for all $a \in S$ and all $\mathfrak p \in \mathcal P^*$, there is a subset $G_S^* $ of $\mathcal C_v (H)$ such that $\mathcal B (G_S^*) \hookrightarrow \mathcal F (G_S^*)$ is a divisor theory and a transfer homomorphism $\theta_S^* \colon S \to \mathcal B (G_S^*)$.

\end{enumerate}
\end{theorem}

\noindent
{\it Remark.} Item 2 is of interest when studying the set $\Delta^* (H)$ of minimal distances of $H$ (see Remark \ref{3.9} and \cite{Ge-Zh16c}).

\begin{proof}
1. By \eqref{v-class-group}, there is a cofinal divisor homomorphism $\mathcal H \to \mathcal I_v^* (H)$, where $\mathcal H = \{aH \mid a \in H\}$ and $\mathcal C_v (H) = \mathsf q (\mathcal I_v^* (H))/\mathsf q (\mathcal H )$. By \cite[Section 5]{Ge-Ka-Re15a}, we have
\[
\mathcal I_v^* (H) \cong D = \mathcal F (\mathcal P) \time T  \,.
\]
If $H_0$ denotes the image of $\mathcal H$ under this isomorphism, then $\mathcal C_v (H) \cong G := \mathsf q (D)/\mathsf q (H_0)$. If $G_{\mathcal P} = \{[p] = p \mathsf q (H_0) \mid p \in \mathcal P\} \subset G$ and $\iota \colon T \to G$, $\iota (t) = t \mathsf q (H_0)$, then Proposition \ref{4.4}.1 implies that  $\boldsymbol \beta_0 \colon H_0 \to B = \mathcal B (G_{\mathcal P}, T, \iota)$ is a transfer homomorphism. The localizations $H_{\mathfrak p}$ are primary and $v$-noetherian for all $\mathfrak p \in \mathcal P^*$. Thus, $F = \mathcal F (G_{\mathcal P}) \times T$ is $v$-noetherian and $B$ is $v$-noetherian by \cite[Proposition 2.4.4]{Ge-HK06a}. Since $\mathsf q (F)/\mathsf q (B) \cong G$ is finite, $B$ is weakly Krull by \cite[Lemma 5.2]{Ge-Ka-Re15a}.

2. Let $S \subset H$ be a divisor-closed submonoid with $\mathsf v_{\mathfrak p} (aH)=0$ for all $a \in S$ and all $\mathfrak p \in \mathcal P^*$. Then $\mathcal S = \{aH \mid a \in S \} \subset \mathcal H$ is a divisor-closed submonoid, and we denote by $S_0 \subset H_0$ the image of $\mathcal S$ under the above isomorphism. There is a subset $P_0 \subset P$ such that $S_0 \hookrightarrow \mathcal F (P_0)$ is a cofinal divisor homomorphism, whence $S_0$ is a Krull monoid whose class group is a subgroup of $G$. Thus, the assertion follows from Proposition \ref{3.3}.
\end{proof}

\medskip
As mentioned in the Introduction, the refinement (done in Proposition \ref{3.3} and in Theorem \ref{T:MainThm}) of the block homomorphism (as given in Lemma \ref{3.2}) simplifies arithmetical investigations. Indeed, this was the original motivation for work in this direction, and it was used a lot for finite abelian groups (see the references given in the Introduction). We discuss here what is known in case of infinite groups. Let $G$ be an abelian group and $G_0 \subset G$ be a subset. If $G_0$ is finite, then $\mathcal B (G_0)$ is finitely generated, which implies a variety of arithmetical finiteness results. In particular, $\mathcal B (G_0)$ is locally tame, has finite catenary degree, and sets of lengths are well-structured. On the other side of the spectrum, we know that if $G_0$ contains an infinite subgroup of $G$, then $\mathcal B (G_0)$ is not locally tame, its catenary degree is infinite, and every finite subset of $\N_{\ge 2}$ occurs as a set of lengths (\cite{Ka99a}). If $G_0$ is infinite but does not contain an infinite subgroup, then information on the arithmetic is available only in very special cases (e.g., \cite{Ha02c, Ge-Gr09b, Ge-Ka10a, De-Ze20a, Gr22a}). The best understood case is when $G$ is an infinite cyclic group. In this case, there is a simple characterization of several arithmetical finiteness properties (including local tameness and the finiteness of the catenary degree) of $\mathcal B (G_0)$ in terms of $G_0$ (\cite[Theorem 4.2]{Ge-Gr-Sc-Sc10}).

We discuss two examples of subsets $G_0 \subset \Z^s$ with $s \ge 2$.  In Example  \ref{4.7} (where $s=2$), we use the strategy developed in Section \ref{3} to obtain a transfer homomorphism to a monoid of zero-sum sequences over a subset of an infinite cyclic group. Thus, in that case Theorem 4.2 of \cite{Ge-Gr-Sc-Sc10} yields information on the arithmetic. In the first example, it will turn out that $\mathcal B (G_0)$ is an inner direct product. If $\boldsymbol 0 \in G_0$, then the sequence $\boldsymbol 0 \in \mathcal A (G_0)$ is a prime element of $\mathcal B (G_0)$, whence $\mathcal B (G_0) = \mathcal B ( \{ \boldsymbol 0\}) \mathcal B (G_0 \setminus \{\boldsymbol 0\})$ is an inner direct product. Thus, when studying the structure of $\mathcal B (G_0)$, we may always assume that $\boldsymbol 0 \notin G_0$.

\smallskip
\begin{example} \label{4.6}
Let $s \ge 2$, $k_2, \ldots, k_s \in \mathbb Q \setminus \{-1\}$, and $G_0 = G_1 \cup G_2 \subset \Z^s$, where
\[
G_1 \subset \{ \boldsymbol a = (a_1, \ldots, a_s) \in \Z^s \mid a_{\nu} = k_{\nu} a_1 \ \text{for all} \ \nu \in [2,s] \}
\]
and
\[
G_2 \subset \{ \boldsymbol b = (b_1, \ldots, b_s) \in \Z^s \setminus \{\boldsymbol 0\} \mid b_2 = -b_1 \} \,.
\]
\begin{enumerate}
\item $G_1 \cap G_2 = \emptyset$ and $\mathcal B (G_0) =  \mathcal B (G_1) \mathcal B (G_2)$ is an inner direct product.

\smallskip
\item $\mathcal B (G_0)$ is locally tame resp. has finite catenary degree if and only if $\mathcal B (G_1)$ and $\mathcal B (G_2)$ are locally tame resp. have finite catenary degree.
\end{enumerate}

\noindent
{\it Remark.} The same equivalence as stated in 2. holds true for several further arithmetical finiteness properties (see, for example, \cite[Theorem 3.11]{Ge-Gr09b}).

\begin{proof}
1.
Let $S \in \mathcal B  (G_0 )$ be given. Then $S = S_1 S_2$, where
\[
S_1 = \prod_{\nu=1}^l \boldsymbol a_{\nu} \in \mathcal F (G_1) \quad \text{and} \quad S_2 = \prod_{\nu=1}^m \boldsymbol b_{\nu} \in \mathcal F (G_2) \,.
\]
It suffices to show that $S_1 \in \mathcal B (G_1)$. Since $S \in \mathcal B (G_0)$, we have
\[
\boldsymbol 0 = \sigma (S) = \sigma (S_1)+\sigma (S_2) = \boldsymbol a_1 + \ldots + \boldsymbol a_l + \boldsymbol b_1 + \ldots + \boldsymbol b_m
\]
and hence
\[
\begin{aligned}
0 & =  a_{1,1} + \ldots + a_{l,1} + b_{1,1} + \ldots + b_{m,1} \quad \text{and} \\
0 & =  a_{1,2} + \ldots + a_{l,2} + b_{1,2} + \ldots + b_{m,2}  \,.
\end{aligned}
\]
Adding the last two equations and using that $b_{j,2}=-b_{j,1}$ for all $j \in [1,m]$, we infer that
\[
0  =  (a_{1,1} + \ldots + a_{l,1}) + (a_{1,2} + \ldots + a_{l,2}) = (a_{1,1} + \ldots + a_{l,1}) + k_2(a_{1,1} + \ldots + a_{l,1})  \,,
\]
and hence $a_{1,1} + \ldots + a_{l,1} = 0$. Thus, for all $\nu \in [2,s]$, we infer that $a_{1, \nu} + \ldots + a_{l, \nu} = k_{\nu}(a_{1,1} + \ldots + a_{l,1}) = 0$, and hence $\boldsymbol 0 = \boldsymbol a_1 + \ldots + \boldsymbol a_l = \sigma (S_1)$.

2. This follows from 1. and from \cite[Proposition 1.6.8]{Ge-HK06a}.
\end{proof}
\end{example}

\smallskip
\begin{example} \label{4.7}
Let  $G_0= G_1 \cup  \{ \boldsymbol a\} \subset \Z^2$, where $G_1\subset\N_0^2$ and $\boldsymbol a=(a_1,a_2)\in\Z^2$, with $a_1<0$ and $a_2<0$, such that $[G_0]=\langle G_0\rangle$ (i.e., $G_0$ is condensed). Since $[G_0]=[G_1]+\N_0  {\boldsymbol a}$
and $\langle G_0\rangle= \langle G_1\rangle +\Z {\boldsymbol a}$ we obtain in particular
\begin{equation}\label{E:sumepi}
[G_1]+\Z {\boldsymbol a}=\langle G_1 \rangle + \Z {\boldsymbol a} \quad .
\end{equation}
We write  $F\in\free(G_0)$ in the form $F={\boldsymbol a}^{\mathsf v_{\boldsymbol a}(F)}F'$ with $F'\in\free(G_1)$. For every $B\in\block(G_0)$, we have
\begin{equation}\label{E:Charblock}
\mathsf v_{\boldsymbol a}(B) \boldsymbol a = -\sigma(B')\quad .
\end{equation}
Hence we get an injective homomorphism
\[
\varphi\colon \block(G_0)\rightarrow \free(G_1),\quad B\mapsto B'\quad .
\]
We set $H=\varphi(\block(G_0))$, whence $\block(G_0)\cong H$, and  $\Gamma =\langle G_1\rangle/(\langle G_1\rangle\cap \Z {\boldsymbol a})$. We define
$\bar{\sigma}\colon \free(G_1)\rightarrow \Gamma$ to be the composition of the sum homomorphism
$\free(G_1)\rightarrow \langle G_1\rangle$ and the canonical homomorphism $\langle
G_1\rangle\rightarrow \Gamma$. It follows from  \eqref{E:sumepi}  that $\bar{\sigma}$ is onto.

We claim that $H=\bar{\sigma}^{-1}(0)$. The inclusion $\subset$ follows from \eqref{E:Charblock}.
Conversely, let $F'\in\bar{\sigma}^{-1}(0)$. Then $\sigma(F')=k \boldsymbol a$ for some $k\in\Z$. Since $\boldsymbol a\in
(-\N)^2$ and $\sigma(F')\in\N_0^2$ we have $k\leq 0$. Then $\boldsymbol a^{-k}F'\in \block(G_0)$ and therefore
$F'\in H$.

To show that $H \hookrightarrow \free(G_1)$ is cofinal, let $F'\in \free(G_1)$ be given. Choose $F''\in
\free(G_1)$ such that $\bar{\sigma}(F'')=-\bar{\sigma}(F')$, whence $F'F''\in H $.

From \cite[Proposition 2.5.1]{Ge-HK06a} we obtain now: $H \hookrightarrow \free(G_1)$ is a cofinal divisor homomorphism and there is a monomorphism $s\colon \free(G_1)/H\rightarrow \Gamma$ such that
\[
\begin{xy}
\Vtriangle/{->}`{->}`{<-}/[{\free(G_1)}`{\Gamma}`{\free(G_1)/H};{\bar{\sigma}}`{\textrm{can}}`{s}]
\end{xy}
\]
is commutative. Since $\bar{\sigma}$ is onto, $s$ is an isomorphism.
Using the block homomorphism associated to the divisor homomorphism $H \hookrightarrow \free(G_1)$ and the
isomorphisms $s$ and $\block(G_0)\cong H$, we obtain a transfer homomorphism
\[
\block(G_0)\rightarrow \block(\bar{\sigma}(G_1)) \,.
\]
We consider two special cases. In both of them, it will turn out that the torsion-free rank of $\Gamma$ equals one.

\smallskip
\noindent
{\bf Special Case 1.} Let $G_1=\N_0^2$ and $\boldsymbol a=(a_1,a_2)$ with $a_1<0$
and $a_2<0$. Since $\N_0^2$ is a G-monoid (see \cite[Chapter 2.7]{Ge-HK06a} for background on G-monoids) and $\N^2$ is the intersection of all non-empty
prime $s$-ideals of $\N_0^2$, we have $\N_0^2+\N_0 {\boldsymbol a}=\Z^2$. With all notation from above, we have that
$\Gamma$ is a finitely generated group with torsion-free rank equal to $1$. Since
$\sigma(\free(\N_0^2))=\N_0^2$ and $\Z^2=\langle \N_0^2\rangle$, we obtain in this case a transfer homomorphism
\[
\block(G_0)\rightarrow\block{(\Gamma)} \,.
\]

\smallskip
\noindent
{\bf Special Case 2.}
Let $G_1$ be set of all $(m,n)\in \N_0^2$ such that $n\geq m^2$. We set
$\boldsymbol a=(-1,-2)$. From
\[
\{ (0,1),(1,1)\}\subset G_1\subset \{ (m,n)\in\N_0^2\st n\geq m\}=[(0,1),(1,1)]
\]
we obtain
\[
[G_1]=[(0,1),(1,1)] \,.
\]
Note that $[(0,1),(1,1)]$ is a factorial G-monoid. Arguing as in example before we obtain
$[G_1]+\N_0 {\boldsymbol a}=\Z^2$, which in particular implies that, $G_0$ is condensed. The homomorphism
\[
\langle G_1\rangle=\Z^2\rightarrow \Z, (u,v)\mapsto 2u-v
\]
has kernel $\Z {\boldsymbol a}$. Hence we obtain an isomorphism $\varphi\colon\Gamma\rightarrow\Z$ such that
$\varphi(\bar{\sigma}(G_1))=\varphi(G_1)$. We show $\varphi(G_1)=-\N_0 \cup \{ 1\}$.  For all $n\in\N_0$,
we have $(0,n)\in\ G_1$ and hence $-n=\varphi(0,n)\in \varphi(G_1)$. It remains to show that
$\varphi(G_1)\cap \N=\{1\}$. Let $(n,m)\in G_1$ be such that $\varphi(m,n)=2m-n>0$. Then $m^2\leq
n<2m$ which implies $m<2$. If $m=0$, then $n<0$, contradiction. If $m=1$ then $n=1$ and
$\varphi(m,n)=1$.
\end{example}

\medskip

\providecommand{\bysame}{\leavevmode\hbox to3em{\hrulefill}\thinspace}
\providecommand{\MR}{\relax\ifhmode\unskip\space\fi MR }
\providecommand{\MRhref}[2]{%
  \href{http://www.ams.org/mathscinet-getitem?mr=#1}{#2}
}
\providecommand{\href}[2]{#2}

\end{document}